\documentclass{amsart}
\usepackage{amssymb}
\usepackage{xspace}
\usepackage{enumerate}
\usepackage[plain]{fancyref}
\usepackage{ifthen}
\usepackage{hyperref}
\usepackage{epsfig}

\def\repeat#1#2 {\expandafter\gdef\csname B#1\endcsname {\mathbb{#1}}
  \ifthenelse{\equal{#2}{*}}{}{\repeat #2 }}
\repeat ABCDEFGHIJKLMNOPQRSTUVWXYZ*
\def\repeat#1#2 {\expandafter\gdef\csname C#1\endcsname {\mathcal{#1}}
  \ifthenelse{\equal{#2}{*}}{}{\repeat #2 }}
\repeat ABCDEFGHIJKLMNOPQRSTUVWXYZ*

\newcounter{last-index}
\newcommand{\xxx}[2][] {%
  \private{\ifthenelse{\equal{#1}{}}
    {\underline{$\bullet$}}{\uline{#1}}\marginpar{\tiny #2}}%
  \public{#1}\xspace}

\newtheorem{thm}{Theorem}[section]
\newtheorem{lem}[thm]{Lemma}
\newtheorem{prop}[thm]{Proposition}

\newtheorem{conj}[thm]{Conjecture}

\theoremstyle{definition}
\newtheorem{defn}[thm]{Definition}

\theoremstyle{remark}
\newtheorem{ex}[thm]{Example}
\newtheorem{rem}[thm]{Remark}
\newtheorem{prob}[thm]{Problem}

\newcommand*\fancyrefthmlabelprefix{thm}\frefformat{plain}{\fancyrefthmlabelprefix}{Theorem~#1}
\newcommand*\fancyreflemlabelprefix{lem}\frefformat{plain}{\fancyreflemlabelprefix}{Lemma~#1}
\newcommand*\fancyrefproplabelprefix{prop}\frefformat{plain}{\fancyrefproplabelprefix}{Proposition~#1}
\newcommand*\fancyrefcorlabelprefix{cor}\frefformat{plain}{\fancyrefcorlabelprefix}{Corollary~#1}
\newcommand*\fancyrefclaimlabelprefix{claim}\frefformat{plain}{\fancyrefclaimlabelprefix}{Claim~#1}
\newcommand*\fancyreffactlabelprefix{fact}\frefformat{plain}{\fancyreffactlabelprefix}{Fact~#1}
\newcommand*\fancyrefconjlabelprefix{conj}\frefformat{plain}{\fancyrefconjlabelprefix}{Conjecture~#1}
\newcommand*\fancyrefdefnlabelprefix{defn}\frefformat{plain}{\fancyrefdefnlabelprefix}{Definition~#1}
\newcommand*\fancyrefquestionlabelprefix{question}\frefformat{plain}{\fancyrefquestionlabelprefix}{Question~#1}
\newcommand*\fancyrefconstlabelprefix{const}\frefformat{plain}{\fancyrefconstlabelprefix}{Construction~#1}
\newcommand*\fancyrefsetuplabelprefix{setup}\frefformat{plain}{\fancyrefsetuplabelprefix}{Setup~#1}
\newcommand*\fancyrefexlabelprefix{ex}\frefformat{plain}{\fancyrefexlabelprefix}{Example~#1}
\newcommand*\fancyrefremlabelprefix{rem}\frefformat{plain}{\fancyrefremlabelprefix}{Remark~#1}
\newcommand*\fancyrefproblabelprefix{prob}\frefformat{plain}{\fancyrefproblabelprefix}{Problem~#1}
\frefformat{plain}{\fancyrefeqlabelprefix}{(#1)}
\newcommand*\fancyrefitemlabelprefix{item}\frefformat{plain}{\fancyrefitemlabelprefix}{(#1)}
\frefformat{plain}{\fancyreffiglabelprefix}{Figure~#1}
\frefformat{plain}{\fancyrefseclabelprefix}{Section~#1}
\newcommand*\fancyrefsubseclabelprefix{subsec}\frefformat{plain}{\fancyrefsubseclabelprefix}{Subsection~#1}
\newcommand*\fancyrefchalabelprefix{cha}\frefformat{plain}{\fancyrefchalabelprefix}{Chapter~#1}

\renewcommand{\proof}{\emph{Proof \/}}
\newcommand{\proofof}[1]{\emph{Proof \/ of #1}}
\newcommand{\proofend}{~\rule{2mm}{3mm}}
\newcommand{\defeq}{\buildrel{\scriptstyle\rm def}\over=\;}
\newcommand{\size}[1]{|#1|}
\newcommand{\group}[2]{\langle#1,#2\rangle}

\newcommand{\result}[1]{\smallskip
\par\noindent{\narrower{\sl#1}\par}\smallskip\par}

\newcommand{\Agroup}{\group{\CA}{\oplus}}
\newcommand{\Reg}{\text{Reg}}

\def\egypont{\egypt}
\def\ketpont{\ketpt}
\def\egypt#1{{\setbox0=\hbox{$\overline{\strut #1}$}%
  \setbox2=\hbox to\wd0{\hss$\bullet$\hss}%
  \dimen0=\ht2 \advance\dimen0 by-.55555pt
    \divide\dimen0 by-2
    \advance\dimen0 by-1.1111pt 
    \advance\dimen0 by\ht0
  \hbox{$\rlap{\raise\dimen0\box2}\overline{\strut #1}$}}}
\def\ketpt#1{{\setbox0=\hbox{$\overline{\strut #1}$}%
  \setbox2=\hbox to\wd0{\hss$\bullet$\hss$\bullet$\hss}%
  \dimen0=\ht2 \advance\dimen0 by-.55555pt
    \divide\dimen0 by-2
    \advance\dimen0 by-1.1111pt 
    \advance\dimen0 by\ht0
  \hbox{$\rlap{\raise\dimen0\box2}\overline{\strut #1}$}}}

\def\harompt#1{{\setbox0=\hbox{$\overline{\strut #1}$}%
  \setbox2=\hbox to\wd0{\hss$\bullet$\hss$\bullet$\hss$\bullet$\hss}%
  \dimen0=\ht2 \advance\dimen0 by-.55555pt
    \divide\dimen0 by-2
    \advance\dimen0 by-1.1111pt 
    \advance\dimen0 by\ht0
  \hbox{$\rlap{\raise\dimen0\box2}\overline{\strut #1}$}}}
\def\negypt#1{{\setbox0=\hbox{$\overline{\strut #1}$}%
  \setbox2=\hbox to\wd0{\hss$\bullet$\hss$\bullet$\hss$\bullet$\hss$\bullet$\hss}%
  \dimen0=\ht2 \advance\dimen0 by-.55555pt
    \divide\dimen0 by-2
    \advance\dimen0 by-1.1111pt 
    \advance\dimen0 by\ht0
  \hbox{$\rlap{\raise\dimen0\box2}\overline{\strut #1}$}}}

\def\EGYIK {{NK81203}}
\def\MASIK {{K84233}}

\begin{document}
\title {On Triple Lines  and Cubic Curves\\ ---  the Orchard Problem revisited}
\author{Gy\"orgy Elekes ~and  Endre Szab\'o}
\address{Gy\"orgy Elekes\\E\"otv\"os University, Budapest}
\email{elekes@cs.elte.hu}
\thanks{Elekes is partially supported by HU-NSF grants
  OTKA T014105 T014302 and T019367}
\address{Endre Szab\'o\\
  R\'enyi Institute of the Hungarian Academy of Sciences}
\email{endre@renyi.hu}
\thanks{Szab\'o is partially supported by
  OTKA grants \EGYIK~and \MASIK}

\begin{abstract}
Planar point sets with many triple lines 
(which contain at least three distinct points of the set) 
have been studied for 180 years, started with  Jackson
\cite{Jackson:1821} and followed by Sylvester \cite{Sylvester:2473}.
Green and Tao~\cite{Green-Tao:sets-with-few-lines:2012}
has shown recently that
the maximum possible number of triple lines for an $n$  element set is
$\lfloor n(n - 3)/6\rfloor + 1$.
Here we address the related problem of describing the structure of the 
\emph{asymptotically near-optimal\/} configurations, i.e., of those 
for which the number of straight lines, which go through three or more points, 
has a quadratic (i.e., best possible) order of magnitude. 
We pose the problem whether such point sets must always be related to cubic
curves.
To support this conjecture we settle various special cases; some of them
(Theorems \ref{thm:NoFourOnAlgebraicThm} and \ref{thm:FourInALineThm})
are also related to the four-in-a-line problem of Erd\H{o}s.
\end{abstract}

\maketitle

\section {Introduction}

Given $n$ point in the plane $\BR^2$,
a line is \emph{$3$-rich}, if it  contains precisely $3$ of the given points. 
One of the oldest problems of combinatorial geometry,
the so-called Orchard Problem,
is to maximise the number of $3$-rich lines
(see  Jackson~\cite{Jackson:1821} and
Sylvester~\cite{Sylvester:2473}).
Sylvester showed that the number of $3$-rich lines is $n^2/6 + \CO(n)$,
and recently
Green and Tao~\cite{Green-Tao:sets-with-few-lines:2012}
have found the precise value of the maximum.

\begin{thm}[Orchard Problem. Green--Tao]
  \label{thm:orch-probl-Green-Tao}
  Suppose that $\CH$ is a finite set of $n$ points in the plane.
  Suppose that $n\ge n_0$ for some sufficiently large absolute
  constant $n_0$. Then there are no more then 
  $\big\lfloor n(n - 3)/6\big\rfloor + 1$
  lines that are $3$-rich, that is they contain precisely $3$ points
  of $\CH$.
\end{thm}

Here we address the related problem of describing the structure of the 
\emph{asymptotically near-optimal\/} configurations, i.e., of those 
for which the number of straight lines, which go through three or more points, 
has a quadratic (i.e., best possible) order of magnitude. 

\begin{defn}
Let $\CH$ be a subset of the plane $\BR^2$.
A straight line $l$ is called a {\sl triple line} with respect to 
$\CH$ if there exist three distinct points 
$P_1, P_2, P_3\in l \cap \CH$.
We shall also  use the notation
$$
\harompt{\CH} \defeq \{l\ ;\ \size{l \cap \CH}\ge 3 \}.
$$
We extended the notion of \emph{triple line},
without any change in the definition,
to subsets of the projective plane.
\end{defn}
Note that $\harompt{\CH}$ is a set of {\sl lines}, not
a set of triples; e.g.\ if $\CH$ is a collinear set of  3 or more  points 
then $|\harompt{\CH}|=1$.

Triple lines are not necessarily $3$-rich (as they may be $4$-rich,
$5$-rich, and so on), hence
\fref{thm:orch-probl-Green-Tao} does not directly bound the size of
$\harompt{\CH}$.
In any case, it is easy to find a (non-sharp) quadric upper bound.
Indeed,  each line  with three points contains three segments 
of the $n \choose 2$ which connect pairs of points of $\CH$,
hence
$$
\big|\harompt{\CH}\big|\le {1 \over 3}{n \choose 2} = {n^2 \over 6}- n/6.
$$
The following examples show four simple configurations for which the 
quadratic order of magnitude can really be attained. 
Two  of them consist of three collinear point sets each, the third one
is located on a conic and a straight line,  while the fourth one on a cubic.
\begin{ex}
\label{ex:FirstExa}
If $\CH_1$, $\CH_2$, $\CH_3$ are three copies of an arithmetic progression on
three equidistant parallel lines then $|\harompt{\CH_1\CH_2\CH_3}|\approx
N^2/18$, where  $N$ denotes the total number of points and 
$\harompt{\CH_1\CH_2\CH_3}$ denotes the set of  lines $l$ such that 
there exist three distinct points $P_i\in{l \cap \CH_i}$ for $i=1,2,3$. \\
(It is slightly better to place a point set of ``double density''
on the middle line.)
\end{ex}
\begin{ex}
\label{ex:SecondExa}
Let  $P_1,P_2,P_3$ be the vertices of a non--degenerate triangle, and 
 $\CH_i$ ($i=1,2,3$) point sets on  the line through
the vertices $P_{i-1}$ and $P_{i+1}$, defined by 
\begin{equation}
\CH_i=\Bigl\{X\ ;\ 
	\frac{\overline{P_{i-1}X}}{\overline{XP_{i+1}}}\in
	\{\pm1,\pm2^{\pm1},\pm4^{\pm1},\ldots,\pm2^{\pm (n-1)}\}\Bigr\},
\label{eq:GEPTriple}
\end{equation}
where $i\pm1$ is used mod 3 in the indices of the $P_i$.
(See \fref{fig:TriangleFig}.)
\begin{figure}[h]
{\epsfig{file=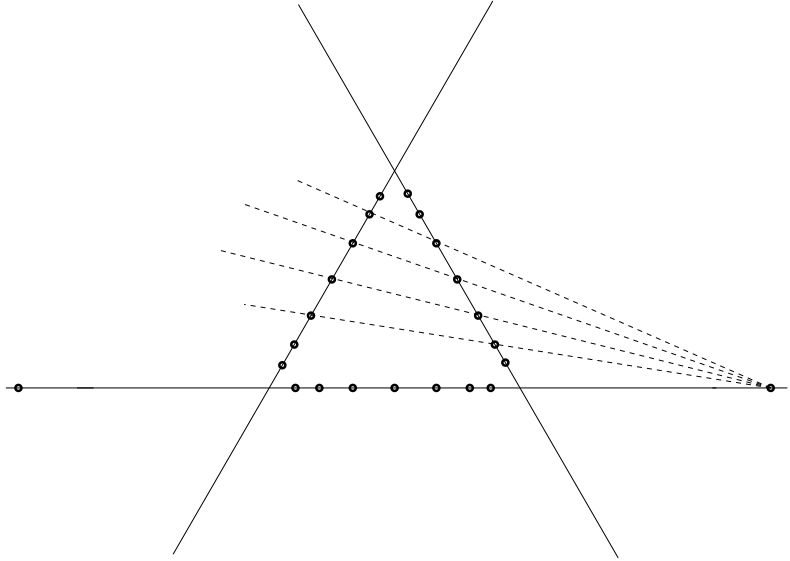}}
\caption{Portion of a triangular configuration with some triple lines marked.}
\label{fig:TriangleFig}
\end{figure}
\par\noindent
Here again $|\harompt{\CH_1\CH_2\CH_3}|\approx
N^2/18$, where  $N$ denotes the total number of points.\\
(The observant reader may have noticed that we allowed $(-1)$ among the ratios,
i.e., $X$ may be a point at infinity.)
\end{ex}
\begin{ex}
\label{ex:ThirdExa}
The $n\choose2$ segments which connect pairs of vertices of a regular $n$--gon
$C$ only determine $n$ distinct slopes. Let $D$ be the set of points on the
line at infinity which correspond to these directions. Then
$|\ketpt{C}\egypont{D}|\approx N^2/8$, where $N=|C\cup D|=2n$ and 
$\ketpt{C}\egypont{D}$ stands for $\harompt{CCD}$.
\end{ex}
\begin{ex}
\label{ex:FourthExa}
The point set $\CH=\{(i,i^3)\ ;\ i=-n,\ldots,n\}$ on the curve $y=x^3$
satisfies $|\harompt{\CH}|\approx N^2/8$, where $N=2n+1$.
This can easily be demonstrated by
making use of the fact that three points $(a,a^3)$, $(b,b^3)$ and $(c,c^3)$
are collinear iff $a+b+c=0$.
\end{ex}

\medskip\par
The goal of this paper
is to show that point sets with many triple lines are,
from several points of view, closely related to cubics.

\section {Problems and results}

\subsection*{A conjecture}

Since all the above examples with a quadratic order of magnitude of the triple
lines involve cubic curves (some of which are degenerate), it is natural to
believe the following. 

\begin{conj}\label{conj:TenConj}
If $|\harompt{\CH}| \ge c\size{\CH}^2 $ 
then ten or more  points of $\CH$ lie on a (possibly degenerate) 
cubic, provided that $\size{\CH}>n_0(c)$.\\ \null \\
\end{conj}
Here the ``magic number'' 10 is the least non-trivial value since any nine
points of $\BR^2$ lie on a cubic. Perhaps even a stronger version may hold:
{\sl 
for every $c>0$ and positive integer $k$ there exist $c^*=c^*(c,k)>0$ and
$n_0=n_0(c,k)$, such that, if
$\;|\harompt{\CH}| \ge c\size{\CH}^2 $
then there is a con-cubic $\CH^*\subset\CH$ with $\size{\CH^*}\ge k$ and 
$|\harompt{\CH^*}| \ge c^*\size{\CH^*}^2 $, 
provided that $\size{\CH}\ge n_0$.
}
\medskip\par\noindent
It is very likely that in place of $k$ above, 
even $c^*\size{\CH}^{\alpha}$ con-cubic points exist
(for some $c^*=c^*(c)>0$ and  $\alpha=\alpha(c)>0$).
An example with only $O(\sqrt{\size{\CH}})$ such points is a 
$k\times k$ square or parallelogram lattice where the points of three parallel
lines provide the set located on a (degenerate) cubic.
Similarly, 
projections of $d$~dimensional cube lattices to $\BR^2$
form structures with only $O(\size{\CH}^{1/d})$ con-cubic points.
\medskip\par
Moreover, if we assume that $\CH$ has no four--in--a--line and
$|\harompt{\CH}| \ge c\size{\CH}^2 $, then perhaps as
many as $c^*\size{\CH}$ of its points will lie on an irreducible cubic.

\subsection*{Results}

In order to support the above conjecture, we 
settle various special cases in the affirmative.
Our main result is the following.
\begin{thm}\label{thm:CubicThm}
In $\BR^2$,
if irreducible algebraic curve of degree $d$
contains a set $\CH$
of $n$ points with $|\harompt{\CH}|\ge cn^2$ then the curve is a cubic ---
provided that $n>n_0(c,d)$. 
\end{thm}
Two simple applications of  the forthcoming slightly more general
\fref{thm:ThreeCurveThm} are the following.

\begin{thm}\label{thm:NoFourOnAlgebraicThm}
In $\BR^2$, no irreducible algebraic curve of degree $d$
can accommodate $n$ points with
$cn^2$ \emph{quadruple} lines if $n>n_0(c,d)$. 
\end{thm}

\begin{thm}\label{thm:FewDirectionsThm}
In $\BR^2$, if a set of $n$ points
located on an irreducible algebraic curve of degree $d$
only determines $Cn$ distinct directions 
then the curve is a conic --- provided that $n>n_0(d,C)$. 
\end{thm}

The above theorems are of algebraic geometric nature,
therefore it is natural to ask analogous questions in complex geometry
(i.e. when the point set and the algebraic curves live in $\BC^2$).
However, in this paper
we restrict our attention to the real plane $\BR^2$.

In some other results (see \fref{sec:ConicSection}) we allow part of the
points (a positive proportion) to be arbitrary and only restrict the rest of
them to a conic.  In this case  it will turn out that
a large subset of the first 
part must be collinear. (Here again, the conic and the straight line,
together, form a degenerate cubic.) The following is the essence of Theorems 
\ref{thm:DegenerateConicThm} and \ref{thm:NonDegenerateConicThm}.
\result{Let $\CH=\CH_1\cup\CH_2$ and assume that $\CH_1$ lies on a (possibly
degenerate) conic $\Gamma$ while $\CH_2\cap\Gamma=\emptyset$.
If $n\le|\CH_1|,|\CH_2|\le Cn$ and $|\ketpt{\CH_1}\egypt{\CH_2}|\ge cn^2$ then 
some $c^*n$ points of $\CH_2$ are collinear. (Here $c^*=c^*(c,C)$ does not
depend on $n$.)}
\noindent
We also mention  a theorem of Jamison  \cite{Ja:84} which can be considered
as another result in the direction of our \fref{conj:TenConj}: 
if the diagonals and sides of a convex $n$--gon only determine 
$n$ distinct slopes (which is smallest possible), then the vertices of the
polygon all lie on an ellipse. In terms of triple lines (and a degenerate
cubic formed by a straight line and an ellipse) this can be formulated
as follows:
\result{(Jamison's Theorem) if\/  $\CH_1$ 
is the vertex set of a convex polygon 
and $\CH_2$ 
lies on the line at infinity 
with $|\CH_1|=|\CH_2|=n$ then
$|\ketpont{\CH_1}\egypont{\CH_2}|={n \choose 2}$
implies that $\CH_1$ lies on an ellipse.
}
\noindent
A similar statement was proven by Wettl \cite{We:87} for finite projective
planes.
\medskip\par\noindent

\subsection*{The structure of the paper}
The aforementioned results
(usually in stronger form) are presented in detail in the last two sections. 
Before that, we list some basic facts on the relation between continuous
curves, collinearity and Abelian groups, concluding in the fundamental
observation \fref{lem:MainLemma}.

\section{Collinearity and groups}\label{sec:CurveSection}

\subsection*{Collinearity on cubics}

\begin{defn}
Let $\Gamma_1$, $\Gamma_2$, $\Gamma_3$ be three (not necessarily
distinct)  Jordan curves (i.e., bijective continuous images of an
interval or a circle)
in the projective plane,
and $\Agroup$ an Abelian topological group.
We say that \emph{collinearity between $\Gamma_1$, $\Gamma_2$ and $\Gamma_3$
can be described by the group operation $\oplus$}, if, for $i=1,2,3$, 
there are homeomorphic monomorphisms (i.e., continuous injections
whose inverses are also continuous)
\[
f_i: \Gamma_i \rightarrow \CA
\]
--- in other words,  ``parametrisation'' of the $\Gamma_i$ with $\CA$ --- 
such that three distinct points $P_1 \in \Gamma_1$, $P_2 \in \Gamma_2$, $P_3 \in \Gamma_3$
are collinear if and only if 
\[
f_1(P_1) \oplus f_2(P_2) \oplus f_3 (P_3) = 0 \in \CA.
\]
\end{defn}
The curves we consider will usually be irreducible components of algebraic
curves in $\BR^2$ --- or subsets thereof.
However, sometimes we must also study general
continuous curves, as well.
\medskip\par
\begin{figure}
{\epsfig{file=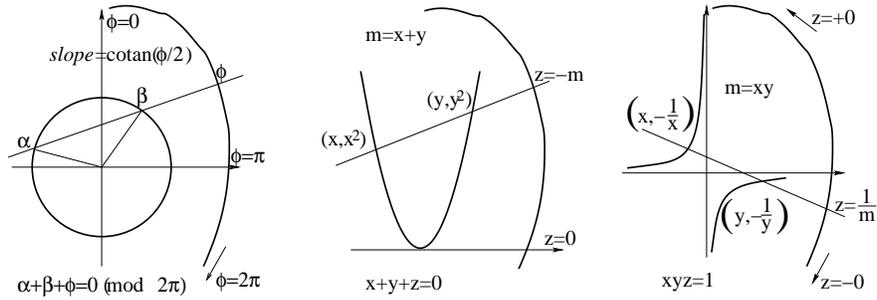}}
\caption{Parametrisation of reducible cubics: a conic plus the line at
	infinity. (Due to lack of sufficient space the line at infinity is
	depicted as a bent curve.)}
\label{fig:ConicFig}
\end{figure}
\begin{figure}
{\epsfig{file=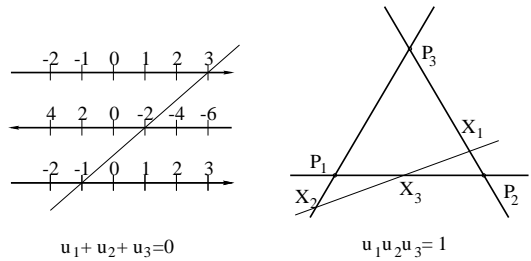}}
\caption{Parametrisation of reducible cubics: three straight lines.
In case of a triangle, 
\ $u_i\defeq f(X_i)=\overline{X_iP_{i-1}}\;/\;\overline{X_iP_{i+1}}$.}
\label{fig:ThreeLinesFig}
\end{figure}

In what follows we denote the set of regular points of an algebraic curve
$\Gamma$ by $\Reg(\Gamma)$. The connected components of $\Reg(\Gamma)$
are Jordan curves.
\begin{prop}\label{prop:CollinearityOnCubics}
  Let $\CC$ be a cubic curve in the projective plane.
  If $\Gamma_1$, $\Gamma_2$, $\Gamma_3$ are (not necessarily distinct)
  connected components of $\Reg(\CC)$,
  then collinearity between them
  can be described by commutative group operation 
  --- unless two of the $\Gamma_i$ are identical straight lines.
\end{prop}
Indeed, for reducible cubics,
Figures~\ref{fig:ConicFig} and \ref{fig:ThreeLinesFig} show
appropriate parametrisation in the real plane. (Any other reducible cubic
is projective equivalent to one of these.)
The groups used are 
\mbox{$\group{\BR}{+}/2\pi\BZ$}, \mbox{$\group{\BR}{+}$},
\mbox{$\group{\BR\setminus\{0\}}{\;\cdot\;}$}
in \fref{fig:ConicFig} 
and \mbox{$\group{\BR}{+}$}, \mbox{$\group{\BR\setminus\{0\}}{\;\cdot\;}$}
in \fref{fig:ThreeLinesFig}, respectively.
If $\Gamma_1=\Gamma_2=\Gamma_3=\CC=\{(x,x^3)\;;\;x\in\BR\}$ then the
parametrisation $f(x,x^3)=x$  works well.
It is also well-known that for irreducible cubics
(i.e. elliptic curves),
suitable parametrisation exist (see, e.g., in \cite{Reid:UAG}). 
\medskip\par

\begin{rem}
Note that in all cases only regular points are parametrised. This will make
no confusion since singular (e.g., multiple) points of a cubic never occur in
proper collinear triples.
\end{rem}

\subsection*{Collinearity on continuous curves}

Throughout this section we consider the graphs of three continuous real
functions. 
\begin{defn}
\label{defn:ssocrfDef}
We call $\alpha$, $\beta$ and $\gamma$ a 
\emph{standard system of continuous real functions} if
\begin{enumerate}[(i)]
\item
they are defined in a neighbourhood $\CD$ of 0;
\item
$\alpha(x)<\beta(x)<\gamma(x)$ for all $x\in\CD$;
\item
any straight line through
any point of the graph of any of the three functions intersects the other two
graphs in at most one point each.
\end{enumerate}
For such functions $\alpha$, $\beta$ and $\gamma$
we denote their graphs (which are Jordan arcs)
by $\overline\alpha$, $\overline\beta$ and $\overline\gamma$.
\end{defn}
\begin{rem}\label{rem:DerivativeRem}
Assumption (iii) is not very strong a requirement;
e.g., if the functions are differentiable at 0 (elsewhere they may not even
be smooth) then $\CD$ can be restricted to a sufficiently small
neighbourhood  of 0 so that (iii) be satisfied there.
\end{rem}
\begin{prop}\label{prop:ContinuityPro}
Let $P(x,\beta(x))$ be a point of the ``middle'' graph $\overline\beta$.
Connect it with lines to the two
points $A_0(0,\alpha(0))$ and $C_0(0,\gamma(0))$; moreover, denote by $C(P)$
and $A(P)$ the points of intersection of these lines with the graphs
$\overline\gamma$ and $\overline\alpha$, respectively (if they exist).
Finally, let $B(P)$ be the intersection of the line through $A(P)$ and $C(P)$
with the graph $\overline\beta$. Then
\begin{enumerate}[(i)]
\item
if $x$ is sufficiently close to 0 then $A(P)$, $B(P)$ and $C(P)$ really exist;
and the composite mappings
\[
x\;\mapsto\; P=P(x,\beta(x)) \; \mapsto \ \begin{cases}
					A(P) \text{ or}\\
					B(P) \text{ or}\\
					C(P) \\
				\end{cases}
\]
are continuous functions $\BR\rightarrow\BR^2$;
\item
for every point $\hat B$ of the graph $\overline\beta$,
sufficiently close to the $y$--axis,  there is a $P$ for which
$\hat B = B(P)$. 
\end{enumerate}
\end{prop}
The straightforward proof using straightforward calculus --- together with the
 Intermediate Value Theorem for (ii) --- is left to the reader.
\proofend
\medskip\par\par\noindent
Next we shall study when will collinearity between
$\overline\alpha$ $\overline\beta$ and $\overline\gamma$
be described by an Abelian topological group $\CA$,
so we will search for parametrisations
$f_\alpha:\overline\alpha\to\CA$,
$f_\beta:\overline\beta\to\CA$ and
$f_\gamma:\overline\gamma\to\CA$.
Part (iii) of \fref{defn:ssocrfDef} also implies that the curves
$\overline\alpha$, $\overline\beta$ and $\overline\gamma$
must be pairwise disjoint. 
That is why, in what follows, we shall only use one notation
$$
f:=(f_\alpha\cup f_\beta\cup f_\gamma)\;:\;
(\overline\alpha\cup\overline\beta\cup\overline\gamma)
\rightarrow A 
$$
in place of three.
\begin{lem}[``Parameter--halving lemma'']\label{lem:ParameterHalvingLemma}
Let $\alpha$, $\beta$ and $\gamma$ form a standard system of continuous real functions. Moreover, let
$B_0=(0,\beta(0))$ and $A_0$, $C_0$, $P$, $A=A(P)$, $B=B(P)$ and $C=C(P)$ 
be as above.
Assume that collinearity between the three graphs is described by 
a group operation
$\group{\CA}{\oplus}$ and mapping (parametrisation) $f$. 
Then
\begin{enumerate}[(i)]
\item
if
\[
\begin{aligned}
f(P)&=f(B_0) \oplus p \text{\quad and }\cr
f(B)&=f(B_0) \oplus b 
\end{aligned}
\]
then $p=b/2$, i.e., $b=p\oplus p$.
\item
if $B$ is sufficiently close to $B_0$ then there really exists a $P$ for which
$f(P)=f(B_0) \oplus b/2$.
\end{enumerate}
\end{lem}
\proof
(i) Note that 
$$
f(A_0) \oplus f(B_0) \oplus  f(C_0) = 0 \in \CA.
$$
Moreover, the collinearity of the triples $C_0PA$ and $CPA_0$
imply
\[
\begin{aligned}
f(A)&=f(A_0)\ominus p;\\
f(C)&=f(C_0)\ominus p,
\end{aligned}
\]
respectively; therefore
\[
\begin{aligned}
f(B)&= p \oplus p \ominus f(A_0)\ominus f(C_0) = \\
	&= p  \oplus p \oplus f(B_0),
\end{aligned}
\]
whence the required identity.
\par
(ii) is obvious from \fref{prop:ContinuityPro}(ii).
\proofend

\subsection*{A fundamental lemma}

The forthcoming \fref{lem:MainLemma}
will work as our first tool for proving
\fref{thm:CubicThm} and the slightly more general 
\fref{thm:ThreeCurveThm}. 
The basic idea is to use the well-known construction
of the group structure on cubics.
If we know a few points on a cubic,
then just by drawing specific lines and marking specific
intersection points we can construct infinitely many new points on that cubic.

The essence of the following statement is that only on
cubics can Abelian groups describe collinearity.

\begin{lem}\label{lem:MainLemma}
  Let $\alpha$, $\beta$, $\gamma$ be a standard system of continuous
  functions defined in a neighbourhood of 0.
  Assume that collinearity between the three graphs is described by 
  a group operation.
  Then their union
  $\overline\alpha\cup\overline\beta\cup\overline\gamma$
  is contained in a (possibly reducible) cubic.
\end{lem}

For the proof we need certain special structures; they will be the topic of
the next subsection.
The proof itself comes then in the subsection afterwards.

\subsection*{Ten point configurations and cantilevers}

Two types of point-line configurations will play special roles in what
follows. The first one consists of ten points and a certain structure of
triple lines while the latter will extend the former one.

Given $\overline\alpha$, $\overline\beta$, $\overline\gamma$
as in \fref{lem:MainLemma}, we define ten
point configurations as follows.

Denote, again, by $A_0$, $B_0$  and $C_0$ the points of intersection of the
$y$--axis with the three graphs, respectively. 

Choose $B_1$ on  $\overline\beta$ sufficiently close to $B_0$ 
in order to make sure that all the forthcoming points exist.
(This will be described later in more detail.) 
Let $A_1$ (resp.~$C_1$) be the point of
intersection of $\overline\alpha$ with the line through $B_1$ and $C_0$
(resp. that of $\overline\gamma$ with the line through $B_1$ and $A_0$).
Define $B_2$ to be the point of intersection of $\overline\beta$ with the
line through $A_1$ and $C_1$.
Let $A_2$ (resp.~$C_2$) be the point of
intersection of   $\overline\alpha$
with the line through $B_2$ and $C_0$
(resp. that of  $\overline\gamma$
with the line through $B_2$ and $A_0$).
\par
\begin{figure}
  {\epsfig{file=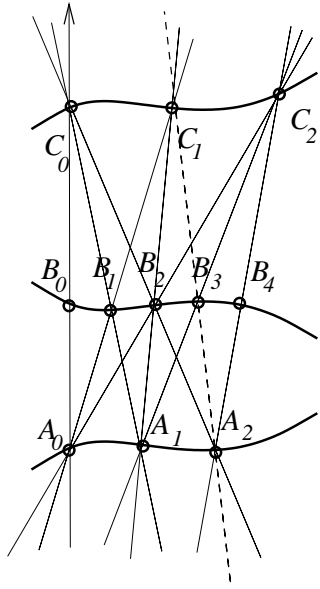}}
  \caption{The straight line $A_2B_3C_1$ is not used in the definition of the
    points.}
  \label{fig:ElevenPoints}
\end{figure}

The definition of $B_3$ is asymmetric: it will be the intersection of
$\overline\beta$
with the line through $A_1$ and $C_2$. Finally, $B_4$ is, again, defined in a
symmetric manner: the intersection of $\overline\beta$
with the line through $A_2$ and $C_2$
(see \fref{fig:ElevenPoints}).
Note that by iterated application of
\fref{prop:ContinuityPro}, the rest
of the points will all exist if $B_1$ is close enough to $B_0$.
\par
The observant reader may have noticed that we defined eleven points
altogether (instead of just ten). 
However, $B_0$ will NOT be in our configuration.
\begin{defn}
\label{defn:TenPtsDef}
Given $\overline\alpha$, $\overline\beta$, $\overline\gamma$
as in \fref{lem:MainLemma}, we call the above
$$\langle A_0,A_1,A_2,B_1,B_2,B_3,B_4,C_0,C_1,C_2\rangle$$
a \emph{ten point configuration}\/ defined by $B_1$.
\end{defn}

\begin{prop}\label{prop:TenPointStructurePro}
If $\alpha$, $\beta$, $\gamma$ 
is a standard system of continuous real functions and
collinearity between their graphs
is described by $\Agroup$ and mapping $f$ then
\begin{enumerate}[(i)]
\item
$A_2$, $B_3$ and $C_1$ are collinear.
\item
More generally, $A_i$, $B_j$ and $C_k$ are collinear iff $i+k=j$.
\item
There is a
$\Delta\in\CA$ such that $f(A_i)=f(A_0)\ominus i\Delta$,
$f(B_i)=f(B_0)\oplus i\Delta$, and 
$f(C_i)=f(C_0)\ominus i\Delta$.
\end{enumerate}
\end{prop}
\proof
Indeed, statement (ii) --- with the exception of (i) --- holds by definition.
For  $\Delta\defeq f(B_1)\ominus f(B_0)$, this implies statement (iii)
by group identities. Finally, (i) follows from (iii), using
 $f(A_0) \oplus f(B_0) \oplus f(C_0)=0$, which, together with (iii), implies
$f(A_2) \oplus f(B_3) \oplus f(C_1)=0$.
\proofend

\begin{lem}[Ten point Lemma]\label{lem:TenPointLemma}
Let $\overline\alpha$, $\overline\beta$, $\overline\gamma$ be a
as in \fref{lem:MainLemma}. Assume, moreover, that a ten point
configuration defined on them is contained in two (possibly reducible) cubics
$\CC_1$ and $\CC_2$. Then $\CC_1=\CC_2$.
\end{lem}
\proof
According to the definition of a standard system of continuous functions,
if a straight line $l$ contains two points of any of the three graphs then $l$
is disjoint from the other two. This leaves us three possibilities for a cubic
$\CC_j$ ($j=1,2$):
\begin{enumerate}[Type 1.]
\item
three straight lines, one through the $A_i$, one through the $B_i$, and one through the $C_i$;
\item
a straight line through all (three or four) points of one of the graphs and a
non-degenerate conic through the rest of them;
\item
an irreducible cubic through all the points.
\end{enumerate}
According to B\'ezout's Theorem \cite{Fulton:AlgCurves}, two distinct
irreducible algebraic curves of degree $k$ and $m$, respectively, can only
intersect in at most $k m$ points. This immediately implies the Lemma.
Indeed, if we assume $\CC_1\ne\CC_2$ for a contradiction, then e.g., if
$\CC_1$ is of type 2 and $\CC_2$  of type 3 then either $\CC_2$ and a
straight line component of $\CC_1$ intersect in four or more points, or
$\CC_2$ and a conic component of $\CC_1$ intersect in seven or more points
--- a contradiction anyway. (The other pairs of types are easier.)
\proofend
\begin{lem}[Nine Point Lemma]\label{lem:NinePointLemma}
Let  $\overline\alpha$, $\overline\beta$, $\overline\gamma$ be a
as in \fref{lem:MainLemma},
consider a ten point configuration on them.
If a (possibly reducible) cubic $\CC$ contains,
with the exception of $B_3$, the
other nine points, then it must also contain $B_3$.
Moreover, all ten points must belong to $\Reg(\CC)$.
\end{lem}
\proof
Define $\delta\defeq f(A_0)\oplus f(B_1) \oplus f(C_0)\in\CA$.
Then $\delta\ne0$ since $A_0$, $B_1$ and $C_0$ are not collinear.
What is $X\in\overline\beta$ for which $f(X)=3\delta$? 
According to \fref{prop:TenPointStructurePro}, it must be the 
point of intersection of the two straight lines $\overline{C_1A_2}$ and 
$\overline{C_2A_1}$.
Finally,
lines passing through a singular point $P\in\CC$, if it has any,
may contain at most two points of $\CC$,
so the lines in our ten point configuration may not pass through $P$.
In particular, $P$ cannot belong to a ten point configuration.
\proofend

\begin{rem}
Note that Lemmas \ref{lem:TenPointLemma} and \ref{lem:NinePointLemma}
also imply 
that two cubics must coincide if they both contain the nine points (with
the exception of $B_3$). However, we shall not need this fact.
\end{rem}

\medskip\par\noindent
Now we extend ten point configurations to what we call
\emph{``cantilevers''}.\footnote
{Cantilever [noun]: a projecting beam or structure supported
only at one end. (The Merriam--Webster Dictionary).}
\\
(We hope that the shape of these structures will really justify this
non-conven\-tional notion.)

Starting from a ten point configuration  on $\alpha$, $\beta$, $\gamma$, we
proceed recursively as follows.

\begin{figure}
  {\epsfig{file=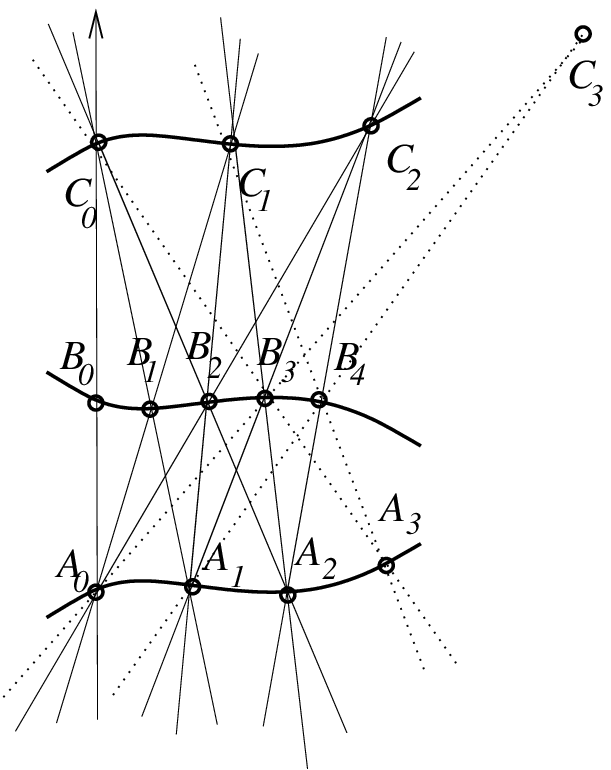}}
  \caption{}
  \label{fig:FourteenPoints}
\end{figure}

Assume that $B_i$ and $B_{i+1}$ have already been defined for an $i\ge 3$. 
Then let $C_i$ be the intersection of the lines $\overline{A_0B_i}$ and 
$\overline{A_1B_{i+1}}$ while 
$A_i$ the intersection of the lines $\overline{C_0B_i}$ and 
$\overline{C_1B_{i+1}}$. Finally, define $B_{i+2}$ to be the 
intersection of $\overline{A_2C_i}$ and $\overline{C_2A_i}$.
(See \fref{fig:FourteenPoints}.)
It is important to note that the construction of cantilevers
use only the ten points, and does not depend on the three curves.
\begin{rem}
Formally, here we work in the projective plane and even allow points of
intersection located on the line at infinity. However, whenever we apply this
construction, all points will lie on the curves 
$\overline\alpha$, $\overline\beta$, and $\overline\gamma$.
\end{rem}

\begin{lem} \label{lem:CantileverOnCurves}
If the straight lines $\overline{A_0B_i}$ and $\overline{A_1B_{i+1}}$
intersect $\overline\gamma$ then this must happen at $C_i$, and similarly
for $\overline{C_0B_i}$, $\overline{C_1B_{i+1}}$, $\overline\alpha$ and $A_i$.
Moreover, if the above intersections all exist (and coincide with 
the $C_i$ and the $A_i$, respectively), then $B_{i+2}$ is located on
$\overline\beta$. 
\end{lem}
\proof
Denote by $X$ and $Y$ the points of intersection of $\overline\gamma$ with
$\overline{A_0B_i}$ and $\overline{A_1B_{i+1}}$, respectively. What is $f(X)$
then? 
By \fref{prop:TenPointStructurePro},
\[
f(X)=\ominus f(A_0) \ominus f(B_i) = 
\ominus f(A_0) \ominus f(B_0) \ominus i\Delta = 
f(C_0)  \ominus i\Delta.
\]
Similarly, $f(Y) = f(C_0) \ominus (i+1-1)\Delta =f(X)$, whence $X=Y$. 
Therefore, also $C_i$ must coincide with these points.\\
A similar argument proves the statement on $B_{i+2}$, too, since in that case
the lines which define it must always intersect $\overline\beta$.
\proofend

\begin{lem} \label{lem:CantileverOnCubics}
  If a cubic $\CC$ contains the nine points
  $A_0$, $A_1$, $A_2$, $B_1$, $B_2$, $B_4$, $C_0$, $C_1$, $C_2$
  of a ten point configuration
  then the entire cantilever (of infinite length)
  built from this configuration
  is contained in $\Reg(\CC)$.
\end{lem}
\proof
By \fref{lem:NinePointLemma} the entire ten point configuration
is contained in $\Reg(\CC)$.
Let $\Gamma_1$, $\Gamma_2$, and $\Gamma_3$
denote the connected components of $\Reg(\CC)$ containing 
$A_0$, $B_1$, and $C_0$, respectively.
By \fref{prop:CollinearityOnCubics} the collinearity between the $\Gamma_i$
is described by a group operation,
let $f_1$, $f_2$, $f_3$ denote the parametrisations.
In this case (i.e. for cubics) all $f_i$ are bijections,
hence they have inverse functions.

Consider the group element $\Delta=f_3(C_1)\ominus f_3(C_0)$.
For all $n\ge 0$ we define the following points on $\CC$:
\begin{eqnarray*}
  A_n' &=& f_1^{-1}\big(f_1(A_0)\ominus n\Delta\big)\\
  B_n' &=& f_2^{-1}\big(f_2(B_1)\oplus (n-1)\Delta\big)\\
  C_n' &=& f_3^{-1}\big(f_3(C_1)\ominus (n-1)\Delta\big)
\end{eqnarray*}
Plugging in $n=0$ and $n=1$ we obtain that
$$
A_0'=A_0,\quad
B_1'=B_1,\quad
C_0'=C_0,\quad
C_1'=C_1.
$$
By assumption $A_0,B_1,C_2$ are collinear, hence 
$f_1(A_0)\oplus f_2(B_1)\oplus f_3(C_2) =0$. This implies that
$$
f_1(A_i')\oplus f_2(B_j')\oplus f_3(C_k') =
\ominus i\Delta\oplus (j-1)\Delta\ominus(k-1)\Delta =
(i+k-j)\Delta
$$
hence $A_i',B_j',C_k'$ are collinear iff $i+k=j$.

Moreover, if a line can intersect $\CC$ in at most three points,
and if two of the intersection points are regular then all of them must
be regular.
Apply this to the line $\overline{C_0B_1}=\overline{C_0'B_1'}$.
The third intersection point of this line with $\Reg(\CC)$
must be $A_1$ by \fref{prop:TenPointStructurePro},
but above we proved it is $A_1'$.
Therefore $A_1'=A_1$.
Similarly, the third intersection point of the line
$\overline{A_1C_1}=\overline{A_1'C_1'}$ with $\Reg(\CC)$
must be $B_2$ on the one hand, and $B_2'$ on the other hand,
which implies $B_2=B_2'$.
Finally apply the same argument to the lines
$\overline{C_0B_2}=\overline{C_0'B_2'}$
and $\overline{A_0B_2}=\overline{A_0'B_2'}$
to obtain that $A_2=A_2'$ and $C_2=C_2'$.

To prove the lemma it is enough to show that $A_n'=A_n$, $B_n'=B_n$
and $C_n'=C_n$ for all $n\ge1$.
We prove it by induction on $n$.
However, it is easier to do the
induction with a slightly stronger statement.
So we shall prove that
$$
A_n'=A_n \;,\quad
B_{n+1}'=B_{n+1} ,\quad
B_{n+2}'=B_{n+2} \;,\quad
C_{n}'=C_{n} 
$$
for all $n\ge 0$.
For $n=0$ we have already seen this.
Assume now that it is true for $n-1$.
Consider the intersection point of the lines
$\overline{C_0B_{n+1}}=\overline{C_0'B_{n+1}'}$
and $\overline{C_1B_{n+2}}=\overline{C_1'B_{n+2}'}$.
On the one hand it must be $A_{n+1}$,
on the other hand it is $A_{n+1}'$,
hence $A_{n+1}'=A_{n+1}$.
Similarly, the intersection point of the lines
$\overline{A_0B_{n+1}}=\overline{A_0'B_{n+1}'}$
and $\overline{A_1B_{n+2}}=\overline{A_1'B_{n+2}'}$
must be $C_{n+1}'=C_{n+1}$.
Finally, the intersection point of
$\overline{C_2A_{n+1}}=\overline{C_2'A_{n+1}'}$
and $\overline{A_2C_{n+1}}=\overline{A_2'C_{n+1}'}$
must be $B_{n+3}=B_{n+3}'$.
This completes the induction step.
\proofend

\subsection*{Proof of \fref{lem:MainLemma}}
\label{subsec:MainLemmaProofSection}

It suffices to show that, for any $x_0$ in the (common) domain $\CD$ of the
functions $\alpha$, $\beta$, and $\gamma$, there exists a cubic $\CC$ which
contains the three graphs \emph{restricted to a sufficiently small
neighbourhood} of $x_0$.
Indeed, if we have such a neighbourhood (for each $x_0$) then it is possible to
extend any of them as follows.
Let $x_1\in\CD$ be one of the endpoints of this neighbourhood (interval) 
and consider a cubic
$\CC_1$ which contains the graphs in a neighbourhood of $x_1$. 
Within the intersection of the two intervals one can find a ten point
configuration contained both by $\CC$ and $\CC_1$. By the Ten Point Lemma
(\fref{lem:TenPointLemma}), $\CC=\CC_1$, i.e., we have a longer neighbourhood
of $x_0$. Thus the \emph{maximal\/} such neighbourhood must be $\CD$
itself.
\par
Now we find an appropriate cubic in a neighbourhood of (without loss of generality) $x_0=0$.
To start with, we select a ten point configuration,  also include $B_0$,
 and extend it to the
other side as follows. Start ``backwards'' from the collinear triple
$A_2$, $B_4$, $C_2$ and define (using $B_3$ in place of the original $B_1$) 
a $5+9+5$ point cantilever --- with $A_0$, $B_0$ and  $C_0$ in the
``middle''. We shall denote this structure by $\CH $.
\par
Define $B_{1/2}$ as in the Parameter Halving Lemma
(\fref{lem:ParameterHalvingLemma}) and, starting from 
$A_0$, $B_0$ and  $C_0$, using this $B_{1/2}$ as reference point,
define  a cantilever with points $A_i$ ($i=0,\ldots,4$),
$B_i$ ($i=0,\ldots,8$) and $C_i$ ($i=0,\ldots,4$).
Of course, the new points will include the old ones, as well, by 
\fref{prop:TenPointStructurePro}(iii).
Also continue the structure ``to the left'' and 
denote this refined (halved) cantilever of 35 points by
$\CH _1$. Keep on defining $B_{1/2^n}$ and $\CH _n$ by  recursive
halving, where 
the latter consists of $(2^{n+2}+1)+(2^{n+3}+1)+(2^{n+2}+1)=2^{n+4}+3$ points.
\par
For each $n$, consider a cubic $\CC_n$ which passes through
$A_0$, $A_{1/2^n}$, $A_{2/2^n}$, $B_{1/2^n}$, $B_{2/2^n}$, $B_{4/2^n}$, $C_0$,
$C_{1/2^n}$, and 
$C_{2/2^n}$. By \fref{lem:CantileverOnCubics}
this cubic contains all points of $\CH_n$.
In particular, all $\CC_n$ must contain the ten point configuration we
started with, hence all these cubics are identical by the Ten Point Lemma
(\fref{lem:TenPointLemma}).
\par
At this point we have a cubic $\CC$ for which 
\[
\bigcup_n\CH _n\subset\CC.
\]
On it, the halving process (starting from $\CH _0$) gives exactly the same
$\CH _n$, whence the parameters which occur in $\cup_n\CH _n$ are
dense somewhere in an open set $\CU$ of the topological group $A$.
Hence so is the point set itself in three corresponding arcs of $\CC$
(i.e., in the homeomorphic pre--images of $\CU$).
By the continuity of $\alpha$, $\beta$, $\gamma$ (and 
$\cup_n\CH_n\subset\CC$), these arcs are completely on $\CC$, as well, 
thus providing the required common
parts. 
\proofend

\subsection*{Surfaces and groups}

Let $F\in \BR[x,y,z]$ be a polynomial of three real variables.
Denote by
\[
S=S_F \defeq \{ (x,y,z)\in\BR^3\ ;\ F(x,y,z)=0 \}
\]
its zero set, i.e., the algebraic surface described by the equation $F=0$.
The degree of $S_F$ is the (total) degree of its defining polynomial $F$.

\begin{defn}
We say that a surface $S\subset\BR^3$ is
\emph{described by a commutative group operation}\/ $\Agroup$
if there are mappings (``parametrisations'') $f_i:\BR\mapsto \CA$
\ for $i=1,2,3$ \ such that
\[
(x_1,x_2,x_3)\in S \ \Leftrightarrow  \ 
f_1(x_1) \oplus f_2(x_2) \oplus f_3(x_3) =0. 
\]
\end{defn}
E.g., the ball of equation $x^2+y^2+z^2= 1$ is described by the additive group
through the mappings $f_i(t)=t^2-1/3$ \ ($i=1,2,3$).

\bigskip\par\par\noindent
One of the main ingredients of our proof is \fref{thm:SpecialSurfaceOrchardThm}
below, proven in \cite{How-to-find-groups}.

  Assume we consider a plane $\alpha x+ \beta y+\gamma z=\delta$, 
intersecting the cube $[0,n]^3$.
  If the coefficients $\alpha, \beta, \gamma, \delta$ are rationals with
  small numerators and denominators then
  this plane will contain $\sim n^2$ lattice points. If we apply
  independent uni-variate
  transformations in the three coordinates, $x,y,z$, then we can
  easily produce 2-dimensional surfaces --- described by some equation
  $f(x)+g(y)+h(z)=\delta$ --- containing a quadratic number of points from a
  product  set $X\times Y \times   Z$, where $|X|=|Y|=|Z|=n$. 
  The main result of \cite{How-to-find-groups} asserts that if
  some appropriate algebraicity conditions hold then (apart from
  being a cylinder) this is the only way for a surface $F(x,y,z)=0$ to contain
  a near--quadratic number of points from such a product set 
  $X\times Y\times Z$.

As usual, we call a function of one or two variable(s)
\emph{analytic\/} at a point if it can be expressed as a convergent power
series in a neighbourhood.
Also, it is analytic on an open set if it is analytic at each of its points.

\begin{thm}[``Surface Theorem'', see \cite{How-to-find-groups}, Theorem~3.]
\label{thm:SpecialSurfaceOrchardThm}

For any positive integer $d$ there exist positive constants
 $\eta=\eta(c,d)$, $\lambda=\lambda(c,d)$  and
 $n_0=n_0(c,d)$ with the following property.\\
If $V\subset\BR^3$ is an algebraic surface (i.e. each component is two
dimensional) of degree $\le d$ then the following are equivalent:
\begin{enumerate}[(a)]
\item \label{item:18}
For at least one $n>n_0(c,d)$ there exist $X,Y,Z \subset\BR$ such that
$|X|=|Y|=|Z|=n$ and 
$$
|V\cap(X\times Y\times Z)|\ge c n^{2-\eta };
$$
\item \label{item:19}
Let $D$ denote the interval $(-1,1)$.
Then either $V$ contains a cylinder over a curve
$F(x,y)=0$ or  $F(x,z)=0$ or  $F(y,z)=0$ or, otherwise,
there are one-to-one analytic functions $f,g,h:D\to\BR$
with analytic inverses
such that $V$ contains the 
$f\times g\times h$-image of a part of the plane
$x+y+z=0$ near the origin:
$$
V\supseteq
\Big\{\,\Big(f(x),g(y),h(z)\Big)\in\BR^3 \ ;\ 
x,y,z\in D\;;\;x+y+z=0\Big\};
$$
\item \label{item:21}
The statement in \eqref{item:19} can be localised as follows.
There is a finite subset $H\subset\BR$ and an irreducible component
$V_0\subseteq V$  such that whenever $P\in V_0$ is a point whose
coordinates are not in $H$, then one may require that
$\Big(f(0),g(0),h(0)\Big)=P$.
\end{enumerate}
\end{thm}
This result indicates a significant ``jump'':
either $V$ has the special form described in \fref{item:19},
in which case a quadratic order of magnitude is possible,
by \fref{item:19}$\Rightarrow$\fref{item:21}; or, else, we cannot
even exceed $n^{2-\eta}$, by \fref{item:18}$\Rightarrow$\fref{item:19}.

\section{Theorems on curves}

Here we present some results on point sets located on algebraic curves and
satisfying certain requirements. 
\par
The first one (\fref{thm:ThreeCurveThm})
is a ``gap version'' of \fref{thm:CubicThm}. It states that there is a
significant difference between cubics and other algebraic curves: on a cubic,
$n$ points can determine as many as $cn^2$ triple lines; otherwise
even as few as $n^{{2-\eta}}$ are impossible for $n$ large enough.
\par
The other result is related to a problem of Erd\H{o}s. He asked if a point set
with $cn^2$ \emph{quadruple\/}  lines must also contain a five-in-a-line. In
\fref{thm:FourInALineThm} we settle this in the affirmative, under the
additional assumption that the points lie on an algebraic curve.
\par
Finally, \fref{thm:FewDirThm} concerns point sets which determine few
distinct directions.

\subsection*{Many triple lines force cubics}

Our first main result states that, of all algebraic curves, only cubics can
accommodate $n$ points with $cn^{2-\eta}$ triple lines.
This is probably far from being best possible; perhaps 
even the existence of as few as $cn^{1+\delta}$ such lines will also imply
the same statement, for any $\delta>0$ and $n>n_0(c,\delta)$.

\begin{thm}\label{thm:ThreeCurveThm}
For every $c>0$ and positive integer $d$ there exist $\eta=\eta(c,d)$ and
 $n_0=n_0(c,d)$ with the following property.
Let $\Gamma_1$, $\Gamma_2$, $\Gamma_3$ be (not necessarily distinct)
irreducible algebraic curves of degree at most $d$ in the plane $\BR^2$.
Assume that $n>n_0$ and 
\begin{enumerate}[(i)]
\item
no two $\Gamma_i$ are identical straight lines;
\item
$\CH_i \subset \Gamma_i$ with $\size{\CH_i} \le n$ ($i=1,2,3$);
\item
 $\size{\egypt{\CH_1}\egypt{\CH_2}\egypt{\CH_3}} \ge cn^{2-\eta}$.
\end{enumerate}
Then  $\Gamma_1\cup \Gamma_2 \cup \Gamma_3$ is a cubic.
\end{thm}

\begin{rem}
If we have an \emph{arbitrary}\/ (i.e., possibly reducible) algebraic curve 
$\Gamma$ of degree $d$ 
and a point set $\CH$ with many triple lines on it, then by the Pigeonhole
Principle, some (at most three) irreducible components of  $\Gamma$ will contain a
subset of $\CH$ which still determines at least $\size{\harompt{\CH}}/d^3$
distinct triple lines. Therefore, the union of these components must be a cubic,
according to the aforementioned Theorem.
\end{rem}
\proofof {\fref{thm:ThreeCurveThm}.}
Let the curves $\Gamma_1$, $\Gamma_2$, $\Gamma_3$ be defined by the polynomial
equations $F_1(x,y)=0$, $F_2(x,y)=0$, $F_3(x,y)=0$, respectively. Three points
$P_i(x_i,y_i)\in\Gamma_i$ ($i=1,2,3$) are collinear iff
\[
F(x_1,y_1,x_2,y_2,x_3,y_3)\defeq 
\begin{vmatrix}
1 & x_1 & y_1 \\
1 & x_2 & y_2 \\
1 & x_3 & y_3 \\
\end{vmatrix}
=0.
\]
Eliminating the $y_i$ from the system of the four equations
\begin{equation}\label{eq:FourEqu}
\bigl\{
F(x_1,y_1,x_2,y_2,x_3,y_3)=0
\bigr\} \cup \bigl\{ 
F_i(x_i,y_i)=0\ \ (i=1,2,3)
\bigr\},
\end{equation}
we get a polynomial relation $f(x_1,x_2,x_3)=0$.
In other words, the projection  to $\BR^3$ (i.e., to the subspace
spanned by the $x_i$ coordinates)
of the two dimensional algebraic variety defined by
\fref{eq:FourEqu} in $\BR^6$, will be contained in the zero-set
of a single polynomial equation $f(x_1,x_2,x_3)=0$.
\par
Let $\eta=\eta(c,d)$ be as in \fref{thm:SpecialSurfaceOrchardThm}.
Denoting the set of the $x$ coordinates of $\CH_i$ by $X_i$ ($i=1,2,3$), we have
that the surface $S_f = \{ f=0 \}$ contains at least $cn^{2-\eta}$ points of 
$X_1\times X_2 \times X_3$.
\par
In other  words, \fref{item:18} of the
Surface Theorem~\ref{thm:SpecialSurfaceOrchardThm} is
satisfied for $V=S_f$ and the $X_i$. Since $S_f$ cannot contain a cylinder by
assumption (i), there exists an irreducible component $V_0\subset S_f$ 
for which also \fref{item:19} --- localised as in \fref{item:21} of the same Theorem --- holds.
\bigskip\par
Pick a generic  point $P(a_1,a_2,a_3)\in V_0\subset S_f$. 
By the definition of the
surface, there exist $b_1$, $b_2$, $b_3\in\BR$ such that, on the one hand,
$Q_i(a_i,b_i)\in\Gamma_i$ for $i=1,2,3$, while on the other hand, these $Q_i$
are collinear. 
We can also assume without loss of generality, that  these three
points are distinct, they are regular points of
$\Gamma_1\cup\Gamma_2\cup\Gamma_3$, and 
the straight line $l$ which contains  them is not tangent to $\Gamma_i$ at
$Q_i$ ($i=1,2,3$).
[Indeed, $V_0$ is two dimensional by
\fref{thm:SpecialSurfaceOrchardThm}\fref{item:19}
while the points to be excluded form a finite number of one dimensional
curves.] 
\par
Moreover, by \fref{item:19} and \fref{item:21} there, collinearity between sufficiently small arcs
of the $\Gamma_i$ around the $Q_i$ is described by $\langle\BR,+\rangle$.
Now if we rotate and/or shift the plane so that $l$ becomes the $y$ axis then,
according to \fref{rem:DerivativeRem}, in a sufficiently small neighbourhood
of 0, the (rotated) $\Gamma_i$ coincide with the graphs of
a standard system of continuous functions.
Thus we can use \fref{lem:MainLemma} to conclude  that a suitable cubic
$\CC$ contains a  non-empty open arc of each $\Gamma_i$. Thus also the union
$\Gamma_1\cup \Gamma_2 \cup \Gamma_3$ of the three \emph{irreducible\/} curves 
is contained in $\CC$.
\par
Finally, they cannot all be contained in a curve of degree $<3$ since in that
case they could not define many triple lines. Therefore, 
$\Gamma_1\cup \Gamma_2 \cup \Gamma_3=\CC$.
\proofend

\subsection*{Four-in-a-line}

Erd\H{o}s \cite{ErPu:IV} posed the problem  whether a set of $n$ points
which contains $cn^2$ collinear four-tuples must also contain five collinear
points. 
To our best knowledge, no progress has been made on this question so far.

In 1995, M.~Simonovits asked the following. \emph{ Is it possible to find $n$
points  on an irreducible algebraic curve of degree 4 which determine $cn^2$
four-in-a-line?\/} (Of course, such a set can contain no five-in-a-line.)
We show here that the answer is in the negative, even in a more general
setting.

\begin{thm}\label{thm:FourInALineThm}
If an algebraic curve $\Gamma$ of degree $d$ accommodates a set $\CH$ of $n$
points with $cn^{2-\eta}$ distinct quadruple lines,
where $\eta=\eta(c,d)$ is the same as in \fref{thm:ThreeCurveThm},
then $\Gamma$ contains four
straight lines, each with $\ge c^{\prime}(c,d)\cdot n^{1-\eta}$ points of $\CH$,
provided that $n>n_0(c,d)$.
\end{thm}

\proof
$\Gamma$ has at most $d$ irreducible components. Classify the $cn^{2-\eta}$
collinear four-tuples (located on  distinct straight lines) according to
which point lies on which component. By the Pigeonhole Principle, some
four (not necessarily distinct) components
$\Gamma_1$, $\Gamma_2$, $\Gamma_3$ and
$\Gamma_4$ generate $cn^{2-\eta}/d^4={c}^{\prime}(c,d)n^{2-\eta}$ quadruple lines.
By \fref{thm:ThreeCurveThm}, any three of the $\Gamma_i$ must form a cubic.
However, this is only possible if they are distinct straight lines.
\proofend

\subsection*{Few directions}

In \cite{EGyRL}, it was shown that if the graph of a polynomial $f\in\BR[x]$
contains $n$ points whose ${n \choose 2}$ connecting lines only determine
a linear number (at most $C n$) distinct directions then the polynomial $f$ is
quadratic.
(Some historic remarks and earlier results concerning sets which determine few
directions can also be found there.) 

Here we extend this to  general algebraic curves.

\begin{thm}\label{thm:FewDirThm}
For every  $C>0$ and positive integer $d$ there is an $n_0=n_0(C,d)$ with the following property.\\
Let $\Gamma_1$ and  $\Gamma_2$ be two (not necessarily distinct) irreducible
algebraic curves, $n>n_0$, and $\CH_i\subset\Gamma_i$ with $\size{\CH_i}=n$
($i=1,2$). 
Assume that among the directions
of  the straight lines $\overline{P_1P_2}$, for $P_i\in\CH_i$ and $P_1\ne
P_2$, 
at most $Cn$ are distinct. Then $\Gamma_1\cup\Gamma_2$ is a
(possibly degenerate) conic.
\end{thm}
\proof
Let $\Gamma_3$ be the line at infinity and $\CH_3$ the set of the 
$\le Cn$ directions on it.
(If someone prefers no points at infinity, they can apply a projective mapping
before proceeding further.)
By assumption, 
$\size{\egypt{\CH_1}\egypt{\CH_2}\egypt{\CH_3}} \ge {n \choose 2}
	>  n^{2-\eta} $ if $n$ is large. 
Hence, by
\fref{thm:ThreeCurveThm}, $\Gamma_1\cup\Gamma_2\cup\Gamma_3$ is a cubic.
Therefore, $\Gamma_1\cup\Gamma_2$ is a conic.
\proofend

\section{Straight lines and conics}\label{sec:ConicSection}

\begin{thm}\label{thm:DegenerateConicThm}
Let $n\le|\CH_1|,|\CH_2|,|\CH_3|\le Cn$ and assume that 
$\CH_1$ and $\CH_2$ lie on the distinct straight lines 
$l_1$ and $l_2$, respectively, while $\CH_3\cap l_1=\CH_3\cap l_2=\emptyset$.
\\
If, moreover, $|\harompt{\CH_1\CH_2\CH_3}|\ge cn^2$, 
then some $c^*n$ of the points of $\;\CH_3$, too,  must be collinear. 
(Here $c^*=c^*(c,C)$ does not depend on $n$.)
\end{thm}

\proof
Apply a projective transform $\pi$ which maps $l_1$ to the line at infinity.
Then some $cn^2$ pairs of points of $\pi(\CH_2)\times\pi(\CH_3)$
determine at most $\size{\pi(\CH_1)}=\size{\CH_1}\le Cn$ distinct directions, 
while $\pi(\CH_2)$ is still collinear.
By a result in \cite{EGy:LinIII} (see Theorem 3 there), also $\pi(\CH_3)$ 
--- hence $\CH_3$, too ---  must contain $c^*n$ collinear points.
~\proofend
\medskip

The following \fref{thm:NonDegenerateConicThm} is the ``elder brother'' of 
\fref{thm:DegenerateConicThm}
in the sense that now we start from a non-degenerate 
conic while the two lines $l_1$, $l_2$ above  can be considered as a
degenerate one. 

\begin{thm}\label{thm:NonDegenerateConicThm}
Let $C>1$ be arbitrary  and $\CH_1$, $\CH_2\subset\BR^2$.
Assume that 
\begin{enumerate}[(a)] 
\item 
$n\le|\CH_1|,|\CH_2|\le Cn$; 
\item 
$\CH_2$ lies on a non-degenerate conic which contains no point of $\CH_1$;
\item $|\egypt{\CH_1}\ketpt{\CH_2}|\ge n^2$.
\end{enumerate}
Then some $c^*n$ of the points of $\CH_1$ must be collinear
(where  $c^*=c^*(C)$ does not depend on $n$.)
\end{thm}
\proof
First, without loss of generality,
we may assume that  every point of $\CH_1$ is incident upon
at least $n$ triple lines. (Otherwise keep on deleting those with less than
$n/(2C)$ such lines and finally, use the new values of 
$n^\prime=n/(2C)$, $C^\prime=2C^2$.)

Moreover, we may assume that the conic which contains $\CH_2$,
is the parabola $y=x^2$. (Else we apply a projective
mapping which maps it to that curve. This can also be done
such a way that no point of $\CH_1$ is mapped to the
line at infinity and the $x$--coordinates of the points in
$\CH_1\cup\CH_2$ become all distinct.)

Denote the coordinates of the points of $\CH_1$ by $(a_i,b_i)$
and the set of the $x$--coordinates of the points of $\CH_2$
by $X$, i.e.,
$$
\begin{aligned}
\CH_1&=\{ (a_i,b_i)\ |\  i=1,2,\ldots,\size{\CH_1}\};\cr
\CH_2&=\{ (x,x^2)\ |\  x\in X\},
\end{aligned}
$$
where, of course,  $|X|=\size{\CH_2}$.

\begin{prop}
Two distinct points $(x,x^2)$,  $(y,y^2)$  of $\CH_2$ and a point 
$(a_i,b_i)\in\CH_1$ are collinear iff
$$
xy-a_ix-a_iy+b_i=0.\proofend
$$
\end{prop}
The above equations can be considered as
functions of type $X\mapsto X$:
$$
y=f_i(x)\defeq\;\frac{a_ix-b_i}{x-a_i}.
$$
These projective mappings $f_i$ are ``vertical projections'' (to $X$)  
of the involutions of the parabola, with centres $(a_i,b_i)$.
\medskip\par
We started with the assumption that  every point of $\CH_1$ is incident
upon at least $n$ triple lines. Therefore, each $f_i$ maps at least $n$
elements of $X$ to elements of $X$. According to \cite{EGyKZ} Theorem 29
(the ``Image Set Theorem''), some $c^*n$ of the $f_i$ must be collinear ---
if we represent them as elements of the three dimensional projective space.
In other words, in that space at least $c^*n$ points of projective coordinates
$(a_i,-b_i,1,-a_i)$ are  combinations of as few as two of them, say
$(a_1,-b_1,1,-a_1)$ and $(a_2,-b_2,1,-a_2)$. Considering the (constant)
third coordinates, this is only possible if --- even as four dimensional
vectors --- 
$(a_i,-b_i,1,-a_i)=\lambda_i(a_1,-b_1,1,-a_1)+(1-\lambda_i)(a_2,-b_2,1,-a_2)$,
for suitable reals $\lambda_i$. We conclude that also the corresponding
$c^*n$ original points 
$P_i(a_i,b_i)\in\CH_1\subset\BR^2$ must be collinear.\proofend

\section{Concluding remarks }
\bigskip\par

Beyond \fref{conj:TenConj} the following remain open.

\begin{prob}
Let $\delta>0$ be arbitrary. Does the conclusion  ``\/$\Gamma_1 \cup \Gamma_2
\cup \Gamma_3 $ is a cubic'' of
\fref{thm:ThreeCurveThm} hold if, in place of (iii), we only assume
$$
\text{(iii*) }\qquad 
\size{\egypt{\CH_1}\egypt{\CH_2}\egypt{\CH_3}} \ge n^{1+\delta}
$$
--- provided that $n>n_0=n_0(\delta,d)$?
\end{prob}

\begin{prob}
Does \fref{thm:FourInALineThm} hold with $n^{1-\eta/2}$ in the statement
(in place of $n^{1-\eta}$)?
\end{prob}

\begin{prob}
Let $\delta>0$ be arbitrary. Does the conclusion  ``\/$\Gamma_1 \cup \Gamma_2$
is a conic''of
\fref{thm:FewDirThm} hold if we only assume that the lines 
$\overline{P_1P_2}$ only determine $\le n^{2-\delta}$ distinct directions
--- in place of $Cn$ --- provided that $n>n_0=n_0(\delta,d)$?
\end{prob}
\bigskip\par

{\bf\qquad Acknowledgements}
\bigskip

We are grateful to Endre Makai for his very constructive comments
on (and simplifications to) some earlier versions of the manuscript
and also to Zolt\'an J\'arai for sharing with us his typesetting 
\TeX{pertise}.  


\end{document}